\documentstyle[11pt]{article}
\date{}
\begin{document}
\setcounter{page}{1}
\renewcommand {\a}{\alpha}
\renewcommand {\b}{\beta}
\newcommand {\g}{\gamma}
\newcommand {\G}{\Gamma}
\newcommand {\h}{\frac{1}{2}}
\renewcommand {\d}{\delta}
\renewcommand {\l}{\lambda}
\renewcommand {\L}{\Lambda}
\renewcommand {\o}{\omega}
\renewcommand {\O}{\Omega}
\newcommand {\ve}{\varepsilon}
\newcommand {\r}{\rho}
\newcommand {\s}{\sigma}
\newcommand {\th}{\theta}
\newcommand {\z}{\zeta}
\newcommand {\vp}{\varphi}
\newcommand {\be}{\begin{equation}}
\newcommand {\ee}{\end{equation}}
\newcommand {\bea}{\begin{eqnarray}}
\newcommand {\eea}{\end{eqnarray}}
\baselineskip 0.3in
\renewcommand{\theequation}{\arabic{section}.\arabic{equation}}
\title{{\bf A remark on a theorem of the Goldbach - Waring type}
\author{ \sc Claus Bauer, Dolby Laboratories, San Francisco, CA\\clausbauer@yahoo.com }
}
\maketitle
\abstract{Let $p_{i}$, $2\leq i\leq 5$ be prime numbers.
It is proved that all but $\ll x^{23027/23040+\epsilon}$ even integers $N\leq x$ can be written as $N=p_{1}^{2}+p_{2}^{3}+p_{3}^{4}+p_{5}^{4}.$}

\section{Introduction and Statement of Results}
\setcounter{equation}{0} In the thirties, I. M. Vinogradov
\cite{v} and Hua \cite{hua} established many fundamental theorems
in additive prime number theory. Their methods were consecutively
applied to various problems in additive number theory. Among
others, Prachar established in
1952, \cite{pr} the following result:\\
{\it There exists a constant \(c>0\) such that all but \(\ll x(\log\,x)^{-c}\) even integers \(N\)
smaller than \(x\) are representable as
\bea
N=p_{1}^{2}+p_{2}^{3}+p_{3}^{4}+p_{5}^{4}
\eea
for prime numbers $p_{i}$.}\\
In \cite{b} and \cite{b1}, this theorem was improved as follows:\\
{\it All but $\ll x^{19193/19200+\epsilon}$ positive even integers
smaller than $x$ can be represented as in (1.1).}\\ Here we
improve
upon this result by showing the following theorem:\\
{\bf Theorem}\,\,\,{\it All but $\ll x^{23027/23040+\epsilon}$
positive even integers
smaller than $x$ can be represented as in (1.1).}\\

\section{\bf Notation and structure of the proof}
\setcounter{equation}{0} We will choose our notation similar as in
\cite{b1}. By \(k\) we will always denote an integer \(k\in
\{2,3,4,5\}\), by \(p\) we denote a prime number and \(L\) denotes
$\log\,x$. \(c\) is a effective positive constant and \(\epsilon\)
will denote an arbitrarily small positive number; both of them may
take different values at different occasions. For example, we may
write
\begin{eqnarray*}
L^{c}L^{c}\ll L^{c},\quad x^{\epsilon}L^{c}\ll
x^{\epsilon}.\end{eqnarray*} $d(n)$ denotes the number of divisors
of $n$ and $[a_{1},..,a_{n}]$ denotes the least common multiple of
the integers $a_{1},..,a_{n}.$ Be further
\[r\sim R\Leftrightarrow R/2<r\leq R,\quad
{\sum\limits_{\chi\,mod\,q}}^{*}=\sum\limits_{\chi\,mod\,q\atop
\chi\,primitive},\quad \quad {\sum\limits_{1\leq a\leq q}}^{*}
=\sum\limits_{1\leq a\leq q\atop (a,q)=1}^{q}.\] We set
\begin{eqnarray*}
P=N^{\frac{13}{180}-\epsilon},\quad Q=NP^{-1}L^{-E}\quad (E>0
\mbox{ will be defined later }),
\end{eqnarray*}
and
\[
\mu=\frac{1}{2}+\frac{1}{3}+\frac{1}{4}+\frac{1}{5}-1.\]
We define for any characters $\chi$, \(\chi_{j}\)\((mod\,q)\), \(q\leq P\)
and a fixed integer $N$:
\[C_{k}(a,\chi)=\sum\limits_{l=1}^{q}\chi(l)
e\left(\frac{al^{k}}{q}\right),\qquad\quad C_{k}(a,\chi_{0})=C_{k}(a,q).\]
\[
Z(q,\chi_{2},\chi_{3},\chi_{4},\chi_{5})= {\sum\limits_{1\leq
h\leq q}}^{*}e\left(\frac{-hN}{q}\right)
\prod\limits_{k=2}^{5}C_{k}(h,\chi_{k}),\]
\[Y(q)=Z(q,\chi_{0},\chi_{0},\chi_{0},\chi_{0}),\quad
A(q)=\frac{Y(q)}{\phi^{4}(q)}
.\] When the variable $N$ is fixed, we will always write $A(q)$ and  neglect the dependency of $A(q)$ on $N$. Otherwise, we will write $A(q,n)$.
\[s(p)=1+\sum\limits_{\alpha\geq 1}A(p^{\alpha}),
\quad S_{k}(\lambda)=\sum_{\sqrt[k]{x}/2^{k+1}
< n\leq \sqrt[k]{x}}\Lambda(n)e(n^{k}{\lambda}),\]
\begin{eqnarray*}S_{k}(\lambda,\chi)=
\sum\limits_{\sqrt[k]{x}/2^{k+1}
\leq n\leq \sqrt[k]{x}}\Lambda(n)\chi(n) e(n^{k}{\lambda}),& &\quad
T_{k}(\lambda)=\sum_{\sqrt[k]{x}/2^{k+1}
\leq n\leq \sqrt[k]{x}} e(n^{k}\lambda),\\
W_{k}(\lambda,\chi)=S_{k}(\lambda,\chi)-E_{0}T_{k}(\lambda,\chi),
\quad & &E_{0}=\left\{\begin{array}{ll}1,& \mbox{if} \,\chi=\chi_{0},\\
0,& \mbox{otherwise}.\end{array}\right\}\end{eqnarray*}
 Using the
circle method we define the major arcs $ M $ and minor arcs $ m $
as follows:
\[ M=\sum\limits_{q\leq P}{\sum\limits_{1\leq a \leq q}}^{*}I(a,q),\mbox{  }I(a,q)=\left[\frac{a}{q}-\frac{1}{Qq},\frac{a}{q}+\frac{1}{Qq}\right], \]
\[m=\left[\frac{1}{Q} ,1+\frac{1}{Q}\right] \setminus M.
\]
Let
\begin{eqnarray*}
R(N)=\sum\limits_{\sqrt[k]{x}/2^{k+1}
\leq n_{k}\leq \sqrt[k]{x},\,k\in\{2,..,5\}\atop
n_{2}^{2}+..+n_{5}^{5}=N}
\Lambda(n_{2})..\Lambda(n_{5}).\end{eqnarray*}
Then we find
\begin{eqnarray}
R(N)=\int\limits_{\frac{1}{Q}}^{1+\frac{1}{Q}}e(-N\alpha)
\prod\limits_{k=2}^{5}S_{k}(\alpha)\,d\alpha & = &
\left(\int\limits_{M}+\int\limits_{m}\right)
e(-N\alpha)\prod\limits_{k=2}^{5}S_{k}(\alpha)\,d\alpha
\nonumber\\
&=:& R_{1}(N)+R_{2}(N).\end{eqnarray} Arguing as in \cite{b1}, we
see that \bea\label{eq:minor} I_{2}(N)\ll N^{\mu}L^{-A}\eea for
any $A>0$ and all but \(\ll
x^{1+2\epsilon}P^{-1/128}<x^{23027/23040+3\epsilon}\) even
integers \(x/2 \leq N < x\). In the sections 3 and 4 we will show
that for any given $A>0$
\begin{eqnarray}
R_{1}(N)=\frac{1}{120}P_{0} \prod\limits_{p\leq P}
s(p)+O\left(x^{\mu} L^{-A}\right),\end{eqnarray} where
\begin{eqnarray}x^{\mu}\ll P_{0}:=\sum\limits_{m_{1}+m_{2}+m_{3}+m_{4}=N\atop
x/2^{k+1}<m_{k}\leq x}\frac{1}{
m^{1-\frac{1}{k}}}
\ll x^{\mu}\quad\mbox{ for }N\in(x/2,x].\end{eqnarray}

Using that
\[\prod\limits_{p\leq P}s(p)\gg (\log P)^{-960},\]
(see p. lemma 4.5 in \cite{b}), the theorem follows from (2.1) -
(2.4).
\section{\bf The major arcs}
\setcounter{equation}{0} We will make use of the following
lemmas:\\
{\bf Lemma 3.1}\,\,\,{\it Let \(f(x),\,g(x)\) and \(f'(x)\) be
three real differentiable and monotonic functions in the interval
\([a,b]\). If  \(\vert f'(x)\vert\leq \theta<1,\,g(x),g'(x)\ll
1,\) then
\[\sum\limits_{a<n\leq b}g(n)e(f(n))=
\int\limits_{a}^{b}g(x)e(f(x))d\,x\,
+O\left(\frac{1}{1-\theta}\right).\]
  }\\{\it Proof:} See lemma 4.8 in \cite{t}.\\
{\bf Lemma 3.2}\,\,\,{\it For primitive characters $\chi_{i}$ mod
$r_{i} $ (i=1,2,3,4) and the principal character $\chi_{0}$ mod
$q$ we have}
\[\sum\limits_{q\leq P\atop r\vert q}\frac{\vert Z(q,\chi_{0}\chi_{1},
\chi_{0}\chi_{2},\chi_{0}\chi_{3},\chi_{0}\chi_{4})\vert}
{\phi^{4}(q)}\,\ll r^{-1+\epsilon}(\log P)^{c},\]
{\it where} $r=[r_{1},r_{2},r_{3},r_{4}].$\\
{\it Proof:}\,\,\,This is lemma 3.3 in \cite{b1}.\\
{\bf Lemma 3.3}
\begin{eqnarray*}
\sum\limits_{q>x}\vert A(n,q)\vert \ll
x^{-1+\epsilon}d(n).\end{eqnarray*} {\it Proof:} The proof follows
literally the proof of lemma of (4.12) in \cite{lz_sq}.\\
{\bf Lemma 3.4}\,\,\,{\it For $P\leq x^{13/80-\epsilon}$ there is
\begin{eqnarray*}\label{eq:sin}
\sum\limits_{N\leq x}\left\vert \prod\limits_{p\leq
P}s(p,N)-\sum\limits_{q\leq P}A(q,N)\right\vert \ll
xP^{-1/3+\epsilon},\end{eqnarray*} which implies that
\begin{eqnarray*}
\prod\limits_{p\leq P} s(p,N)=\sum\limits_{q\leq P}
A(q,n)+0(x^{-\epsilon})\end{eqnarray*} for all but $\ll
x^{1+2\epsilon}P^{-1/3}$ even integers $N$ with} $1\leq N\leq
x.$\\{\it Proof:} This theorem is stated in \cite{b1} for all
$P\leq x^{7/150-\epsilon}$. The proof shows however that it holds
for $P\leq x^{13/80-\epsilon}$ as well.\\
\\Splitting the summation over \(n\) in residue classes modulo
\(q\), we obtain
\[S_{k}\left(\frac{a}{q}+\lambda\right)
=\frac{C_{k}(a,q)}{\phi(q)}T_{k}(\lambda)
+\frac{1}{\phi(q)}\sum\limits_{\chi\,mod\,q}
C_{k}(a,\chi)W_{k}(\lambda,\chi)+O(L^{2}).\]
Thus we obtain from (2.1)
\begin{eqnarray}R_{1}(N)=R_{1}^{m}(N)+R_{1}^{e}(N)+O(x^{\mu} L^{-A})
\quad (\mbox{for any } G>0),\end{eqnarray}
where
\begin{eqnarray*}
R_{1}^{m}(N)&=& \sum\limits_{q\leq P}\frac{1}{\phi^{4}(q)}
{\sum\limits_{1\leq a\leq q}}^{*} \int\limits_{-1/Qq}^{1/Qq}
\prod\limits_{k=2}^{5}C_{k}(a,q)e\left(-\frac{a}{q}N\right)
T_{k}(\lambda)e(-\lambda N)\, d\lambda,\end{eqnarray*}
\begin{eqnarray*}
& & R_{1}^{e}(N)\\
&=&\sum\limits_{k=2}^{5}\sum\limits_{q\leq P}\frac{1}{\phi^{4}(q)}
{\sum\limits_{1\leq a\leq q}}^{*} \int\limits_{-1/Qq}^{1/Qq}
\prod\limits_{l=2\atop l\ne k}^{5} C_{l}(a,q)T_{l}(\lambda)
\sum\limits_{\chi\,mod\,q} C_{k}(a,q)W_{k}(\lambda,\chi)
e\left(-\frac{a}{q}N-\lambda N\right)d\lambda\\
&+&\sum\limits_{k,l=2\atop k<l}^{5} \sum\limits_{q\leq
P}\frac{1}{\phi^{4}(q)} {\sum\limits_{1\leq a\leq q}}^{*}
\int\limits_{-1/Qq}^{1/Qq} \prod\limits_{m\in \{k,l\}}
C_{m}(a,q)T_{m}(\lambda) \prod\limits_{{o=2\atop o\ne k}\atop o\ne
l}^{5}\\
&\times
&\sum\limits_{\chi\,mod\,q}C_{o}(a,\chi)W_{o}(\lambda,\chi)
e\left(-\frac{a}{q}N-\lambda N\right)d\lambda\\
&+&\sum\limits_{k=2}^{5} \sum\limits_{q\leq
P}\frac{1}{\phi^{4}(q)} {\sum\limits_{1\leq a\leq q}}^{*}
\int\limits_{-1/Qq}^{1/Qq} C_{k}(a,q)T_{k}(\lambda)
\prod\limits_{l=2\atop l\ne k}^{5} \sum\limits_{\chi\,mod\,q}
C_{l}(a,q)W_{l}(\lambda,\chi)
e\left(-\frac{a}{q}N-\lambda N\right)d\lambda\\
&+&\sum\limits_{q\leq P}\frac{1}{\phi^{4}(q)} {\sum\limits_{1\leq
a\leq q}}^{*} \int\limits_{-1/Qq}^{1/Qq}
\prod\limits_{k=2}^{5}\sum\limits_{\chi\,mod\,q}
C_{k}(a,\chi)W_{k}(\chi,\lambda)e\left(-\frac{a}{q}N-\lambda
N\right)\,
d\lambda,\\
&=:& S_{1}+S_{2}+S_{3}+S_{4}.\end{eqnarray*} We first calculate
$R_{1}^{m}(N)$. Applying lemma 3.1 yields
\begin{eqnarray*}
T_{k}(\lambda) &=&
\int\limits_{\sqrt[k]{x}/2^{k+1}}^{\sqrt[k]{x}}e(\lambda u^{k})du\,+\,O(1)
=\frac{1}{k}\int\limits_{x/2^{k+1}}^{x}v^{\frac{1}{k}-1}e(\lambda v)
dv\,+O(1)\\
&=& \frac{1}{k}
\sum\limits_{x/2^{k+1}<m\leq x}\frac{e(\lambda m)}{m^{1-\frac{1}{k}}}\,
+O(1).\end{eqnarray*}
Substituting this in \(R_{1}^{m}(N)\) we see
\begin{eqnarray*}
R_{1}^{m}(N)&=&
\frac{1}{120}\sum\limits_{q\leq P}A(q)
\int\limits_{-1/Qq}^{1/Qq}\prod\limits_{k=2}^{5}\left(
\sum\limits_{x/2^{k+1}<m\leq x}\frac{e(\lambda m)}{m^{1-\frac{1}{k}}}\right)
e(-N\lambda)d\lambda\\
&+& O\left(\max_{2\leq l\leq 5}\sum\limits_{q\leq P}\left\vert A(q)
\int\limits_{1/Qq}^{-1/Qq}\prod\limits_{k=2\atop k\ne l}^{5}
\sum\limits_{x/2^{k+1}<m\leq x}\frac{e(\lambda m)}{m^{1-\frac{1}{k}}}
d\lambda\right\vert\right).\end{eqnarray*}
Using lemma 3.3 and the trivial bound
\begin{eqnarray}
\sum\limits_{x/2^{k+1}<m\leq x}\frac{e(\lambda
m)}{m^{1-\frac{1}{k}}} \ll \min\left(
\sqrt[k]{x},\frac{1}{x^{1-\frac{1}{k}} \vert
\lambda\vert}\right)\,\, ,
\end{eqnarray}
we derive using lemma 3.4,
\begin{eqnarray}
R_{1}^{m}(N)&=&
\frac{1}{120}\sum\limits_{q\leq P}A(q)
\int\limits_{-1/2}^{1/2}\prod\limits_{k=2}^{5}\left(
\sum\limits_{x/2^{k+1}<m\leq x}\frac{e(\lambda m)}{m^{1-\frac{1}{k}}}\right)
e(-N\lambda)d\lambda\nonumber\\
&+& O\left(\sum\limits_{q\leq P}
\left\vert A(q)\right\vert
\int\limits_{1/Qq}^{1/2}\frac{1}{x^{3-\mu}\vert \lambda\vert^{4}}d\lambda
\right)
+O(x^{\mu}L^{-A})\nonumber\\
&=&\frac{1}{120}P_{0}
\sum\limits_{q\leq P}A(q)
+O((PQ)^{3}x^{\mu -3}L^{c})+O(x^{\mu}L^{-A})\nonumber\\
&=&\frac{1}{120}P_{0} \prod\limits_{p\geq 1}s(p)
+O(x^{\mu}L^{-A}),\end{eqnarray} for all but $x^{1+2\epsilon}P^{-1/3}$ integers $N\leq x$, where \(P_{0}\)
is defined as in (2.4) and \(E\) is chosen sufficiently large in $
Q=NP^{-1}L^{-E}.$ In the sequel $ E=E(G)$ is fixed. Now we
estimate the terms $S_{i},\,i=1,2,3,4.$ Using lemma 3.3 we can
estimate \(S_{4}\) in the following way:
\begin{eqnarray*}
& &\left\vert S_{4}\right\vert \nonumber\\
&=&\left\vert\sum\limits_{q\leq P}
\frac{1}{\phi^{4}(q)}
\sum\limits_{\chi_{2}\,mod\,q}
\sum\limits_{\chi_{3}\,mod\,q}
\sum\limits_{\chi_{4}\,mod\,q}
\sum\limits_{\chi_{5}\,mod\,q}
Z(q,\chi_{2},\chi_{3},\chi_{4},\chi_{5})
\int\limits_{-1/Qq}^{1/Qq}\prod\limits_{k=2}^{5}
W_{k}(\l,\chi_{j})e(-n\l)d
\,\l\right\vert\\
&\leq & \sum\limits_{r_{2}\leq P} \sum\limits_{r_{3}\leq P}
\sum\limits_{r_{4}\leq P} \sum\limits_{r_{5}\leq P\atop
[r_{2},r_{3},r_{4},r_{5}]\leq P}
{\sum\limits_{\chi_{2}\,mod\,r_{3}}}^{*}{\sum
\limits_{\chi_{3}\,mod\,r_{3}}}^{*}
{\sum\limits_{\chi_{4}\,mod\,r_{4}}}^{*}
{\sum\limits_{\chi_{5}\,mod\,r_{5}}}^{*}\nonumber\\
&\times
&\int\limits_{-1/Q[r_{2},r_{3},r_{4},r_{5}]}^{1/Q[r_{2},r_{3},r_{4},r_{5}]}
\prod\limits_{k=2}^{5}\left\vert W_{k}(\l,\chi_{k})\right\vert
d\,\l \sum\limits_{q\leq P\atop [r_{2},r_{3},r_{4},r_{5}]\vert q}
\frac{\left\vert Z(q,\chi_{2}\chi_{0},\chi_{3}
\chi_{0},\chi_{4}\chi_{0},\chi_{5}\chi_{0})\right\vert}{\phi^{4}(q)}
,\\
&\ll&  L^{c}
 \sum\limits_{r_{2}\leq P}
\sum\limits_{r_{3}\leq P} \sum\limits_{r_{4}\leq P}
\sum\limits_{r_{5}\leq P} [r_{2},r_{3},r_{4},r_{5}]^{-1+\epsilon}
{\sum\limits_{\chi_{2}\,mod\,r_{2}}}^{*}
{\sum\limits_{\chi_{3}\,mod\,r_{3}}}^{*}
{\sum\limits_{\chi_{4}\,mod\,r_{4}}}^{*}
{\sum\limits_{\chi_{5}\,mod\,r_{5}}}^{*}\\
&\times &
\int\limits_{-1/Q[r_{2},r_{3},r_{4},r_{5}]}^{1/Q[r_{2},r_{3},r_{4},r_{5}]}
\prod\limits_{k=2}^{5}\vert W_{k}(\l,\chi_{k})\vert
d\,\l.\end{eqnarray*} Using \([r_{2},r_{3},r_{4},r_{5}] \geq
r_{2}^{6\epsilon}r_{3}^{1/13-2\epsilon}r_{4}^{4/13-2\epsilon}r_{5}^{8/13-2\epsilon}\),
 we obtain
\begin{eqnarray} & & S_{4}\nonumber\\
&\ll& L^{c} \sum\limits_{r_{2}\leq
P}r_{2}^{-\epsilon}{\sum\limits_{\chi_{k} \,mod\,r_{k}}}^{*} \max
_{\vert \l\vert \leq 1/r_{2}Q}\vert W_{2}(\l,\chi_{2}\vert
\sum\limits_{r_{3}\leq P} r_{3}^{-1/13+2\epsilon}
{\sum\limits_{\chi_{3}\,mod\,r_{3}}}^{*}\max _{\vert \l\vert \leq
1/r_{3}Q} \vert W_{3}(\l,\chi_{l}\vert \nonumber\\ &\times &
\sum\limits_{r_{4}\leq
P}r_{4}^{-4/13+2\epsilon}{\sum\limits_{\chi_{4}\,mod
\,r_{4}}}^{*}\left (\int\limits_{-1/Qr_{4}}^{1/Qr_{4}}\vert
W_{4}(\l,\chi_{4}\vert^{2}d\,\l
\right)^{1/2}\nonumber\\
&\times&\sum \limits_{r_{5}\leq
P}r_{5}^{-8/13+2\epsilon}{\sum\limits_{\chi_{5}\,mod \,r_{5}}}^{*}
\left(\int\limits_{-1/Qr_{5}} ^{1/Qr_{5}}\vert
W_{5}(\l,\chi_{5}\vert^{2}d\,\l\right)^{1/2}
\nonumber\\
&\ll & L^{c}  I_{2} I_{3} W_{4}W_{5},\end{eqnarray} where
\begin{eqnarray*}I_{k}&=&\sum\limits_{r\leq P}
r^{-a_{k}} {\sum\limits_{\chi}}^{*}\max_{\vert\l\vert\leq 1/rQ}\vert W_{k}(\l,\chi\vert,\\
W_{k}&=&\sum\limits_{r\leq P}r^{-a_{k}}{\sum\limits_{\chi}}^{*}
\left(\int\limits_{-1/Qr}^{1/Qr}\vert
W_{k}(\l,\chi\vert^{2}d\,\l\right) ^{1/2},\\
a_{k}&=&\left\{\begin{array}{ll}\epsilon,& \mbox{for} \,k=2,\\
\frac{1}{13}-2\epsilon,& \mbox{for}\,k=3,\\
\frac{4}{13}-2\epsilon,& \mbox{for}\,k=4,\\
\frac{8}{13}-2\epsilon,& \mbox{for}\,k=5.\end{array}\right\}
\end{eqnarray*}
Arguing similarly we obtain \bea S_{1}+S_{2}+S_{3} &\ll &
L^{c}\max_{2\leq k,l,m,n\leq 5\atop k+l+m+n=\mu+1}
\max\limits_{\vert\l\vert\leq 1/Q}\left\vert T_{k}(\l)\right\vert
\max\limits_{\vert\l\vert\leq 1/Q}\left\vert T_{l}(\l)\right\vert
\left(\int\limits_{-1/Q}^{1/Q}\left \vert
T_{m}(\l)\right\vert^{2}\,d\l\right)^{1/2}
W_{n}\nonumber\\
&+&L^{c}\max_{2\leq k,l,m,n\leq 5\atop k+l+m+n=\mu+1}
\max\limits_{\vert\l\vert\leq 1/Q}\left\vert T_{k}(\l)\right\vert
\max\limits_{\vert\l\vert\leq 1/Q}\left\vert T_{l}(\l)\right\vert
W_{m}W_{n}\nonumber\\
&+& L^{c}\max_{2\leq k,l,m,n\leq 5\atop k+l+m+n=\mu+1}
\max\limits_{\vert\l\vert\leq 1/Q}\left\vert T_{k}(\l)\right\vert
I_{l} W_{m}W_{n}. \eea We have trivially
\[\max_{\vert\l\vert\leq 1/Q}\left\vert T_{k}(\l)\right\vert\ll x^{1/k}.\]
Using (3.2) we obtain
\[\left(\int\limits_{-1/Q}^{1/Q}\left
\vert T(\l)\right\vert^{2}\,d\l\right)^{1/2}\ll
x^{\frac{1}{k}-\frac{1}{2}}.\] Thus we see from (3.1) and (3.3) -
(3.5) that the proof of (2.3) reduces to the proof of the
following two lemmas:\\{\bf Lemma 3.5}\,\,\,{\it If $P\leq
x^{\frac{13}{180}-\epsilon}$ and $2\leq k\leq 5$}
\[W_{k}\ll_{B} x^{1/k -1/2}L^{-B}\]{\it for any} \(B>0.\)\\
\\{\bf Lemma 3.6}\,\,\,{\it If
$P\leq x^{\frac{13}{180}-\epsilon}$ and $2\leq k\leq 5$}
\[I_{k}\ll x^{1/k}L^{A}\]{\it for a certain
\(A>0\).}\\

\section{Proof of lemma 3.5}
\setcounter{equation}{0}

In order to prove the lemma it is enough to show that \bea
W_{k,R}\ll
x^{\frac{1}{k}-\frac{1}{2}}R^{a_{k}}L^{-B},\eea where
\[W_{k,R}=\sum\limits_{r\sim R }{\sum\limits_{\chi}}^{*}
\left(\int\limits_{-1/Qr}^{1/Qr}\vert W_{k}(\l,\chi\vert^{2}d\,\l\right)
^{1/2}\]
for $R\leq P/2$.
Applying lemma 1, \cite{g} we see
\begin{eqnarray}
& &\int\limits_{-1/Qr}^{1/Qr}\vert W_{k}(\l,\chi)\vert^{2}\,d\,\l\nonumber\\
&\ll & (QR)^{-2}\int \limits_{x/2^{k+2}}^{x}
\left\vert\sum\limits_{t<m^{k}\leq t+Qr \atop x/2^{k+1}<m^{k}\leq
x}\L(m)\chi(m)-E_{0} \sum\limits_{t<m^{k}\leq t+Qr \atop
x/2^{k+1}<m^{k}\leq x}1\right\vert^{2}dt.\end{eqnarray} We set
$X=\max(x/2^{k+1},t)$ and $X+Y=\min(x,t+Qr).$ In the sequel we
will treat the cases $R>L^{D}$ and $ R\leq L^{D}$ for a
sufficiently large constant $D>0$ separately. In the first case we
apply a slight modification of Heath-Brown's identity (\cite{hb})
\begin{eqnarray*}-\frac{\z`}{\z}(s)
&=&\sum_{j=1}^{K} {K \choose j} (-1)^{j-1}\z`(s)
\z^{j-1}(s)M^{j}(s)-\frac{\z`}{\z}(s)
(1-\z(s)M(s))^{K},\end{eqnarray*}with $K=5$ and
\[M(s)=\sum\limits_{n\leq x^{1/5k}}
\mu(n)\]
to the sum
\[\sum\limits_{X<m^{k}\leq X+Y}.\]
Arguing exactly as in part III, \cite{z} we find by applying
Heath Brown's identity and
Perron's summation formula (see \cite{t}, Lemma 3.12) that the
inner sum of (4.3) - where always $E_{0}=0$ because of
$R>L^{D}$ and the primitivity of the characters - is a linear
combination of \(O(L^{c})\) terms of the form
\[S_{k,I_{a_{1}},..,I_{a_{10}}}=\frac{1}{2\pi i}
\int\limits_{-T}^{T}
F_{k}(\frac{1}{2}+iu,\chi)\frac{(X+Y)^{\frac{1}{k}(\frac{1}{2}+iu)}-
X^{\frac{1}{k}(\frac{1}{2}+iu)}}
{\frac{1}{2}+iu}d\,u
+O(T^{-1}x^{\frac{1}{k}+\epsilon}),\]
where \(2\leq T\leq x\),
\[F_{k}(s,\chi)=\prod\limits_{j=1}^{10}f_{k,j}
(s,\chi),\quad f_{k,j}(s,\chi)=\sum\limits_{n\in I_{k,j}}a_{k,j}(n)\chi_{n}
n^{-s},\]
\[ a_{k,j}(n)=\left\{\begin{array}{ll} \log\,n\mbox{ or } 1,&
j=1,\\
1,& 1<j\leq 5,\\
\mu(n),& 6\leq 10.\end{array}\right\},
\quad
I_{j}=(N_{k,j},\,2N_{k,j}],\quad 1\leq j\leq 10,\]
\bea\sqrt[k]{x}\ll
\prod\limits_{j=1}^{10}N_{k,j}\ll \sqrt[k]{x},\quad
N_{k,j}\leq x^{1/5k},\quad 6\leq j\leq 10.\eea

Since
\[
\frac{(X+Y)^{\frac{1}{k}(\frac{1}{2}+iu)}-
X^{\frac{1}{k}(\frac{1}{2}+iu)}}
{\frac{1}{2}+iu}\ll
 \min\left(QRx^{\frac{1}{2k}-1},x^{\frac{1}{2k}}\left(\vert u\vert
+1\right)^{-1}\right)\]
by taking  \(T=x^{2\epsilon}P^{2}(1+\vert\lambda\vert x)\)
and \(T_{0}=x(QR)^{-1},\)
we conclude that
$S_{I_{a_{1}},..,I_{a_{11}}}$
is bounded by
\begin{eqnarray*}
&\ll &
QRx^{\frac{1}{2k}-1}
\int\limits_{-T_{0}}^{T_{0}}
\left\vert F_{k}(\frac{1}{2}+it,\chi)\right\vert
d\,u
+x^{\frac{1}{2k}}
\int\limits_{T_{0}\leq \vert u\vert\leq T}
\left\vert F_{k}(\frac{1}{2}+it,\chi)\right\vert
\frac{d\,u}{\vert u\vert}\\
&+& x^{\frac{1}{k}}P^{-2},\end{eqnarray*} Thus we derive from
(4.2) that in order to prove (4.1) it is enough to show that
\bea\label{eq:45} \sum\limits_{r\sim
R}{\sum\limits_{\chi}}^{*}\int \limits_{0}^{T_{0}}\left\vert
F_{k}(\frac{1}{2}+it,\chi)\right\vert\,dt &\ll&
x^{1/2k}R^{a_{k}-\epsilon}L^{-B},\\
\label{eq:46}\sum \limits_{r\sim R}{\sum\limits_{\chi}}^{*}
\int\limits_{T_{1}}^{2T_{1}} \left\vert
F_{k}(\frac{1}{2}+it,\chi)\right\vert \,dt &\ll& x^{1/2k-1}Q
R^{1+a_{k}-\epsilon}T_{1}L^{-B},\, T_{0}<\vert T_{1}\vert\leq T.
\eea  For the proof of (\ref{eq:45}) and (\ref{eq:46}) we will
prove two propositions. We will need the estimate \bea
\sum\limits_{n\leq x} d^{k}(n)\ll_{k} xL^{c(k)} .\eea We now
establish\\{\bf Proposition \(1\)} \,\,\, {\it If there exists
$N_{k,j_{1}}$ and $N_{k,j2}$ $(1\leq j_{1},j_{2}\leq 5)$ such that
$N_{k,j1}N_{k,j2}\geq P^{2-2a_{k} +3\epsilon}$ then (4.4) is
true.}
\\
{\it Proof:}\,\,\, We suppose without loss of generality
\(j_{1}=1,\,a_{1}(n)=\log\,n\) and $j_{2}=2,\,a_{2}(n)=1.$ Arguing
exactly as in the proof of proposition 1 in \cite{z}, we find
\[f_{k,1}\left(\frac{1}
{2}+it,\chi\right)\ll L\left(\int\limits_{-x^{1/k}}^{x^{1/k}}\left
\vert
L'\left(\frac{1}{2}+it+iv,\chi\right)\right\vert^{4}\frac{dv}{1+
\vert v\vert}\right)^{1/4}+L,\]and so we find by using lemma 3.7:
\begin{eqnarray*}
& &\sum\limits_{r\sim R}{\sum\limits_{\chi}}^{*}
\int\limits_{0}^{T_{0}}\left\vert f_{1}\left(\frac{1}{2}+it,\chi
\right)\right\vert^{4}\,dt\\&\ll &    L^{4}\int\limits_{-x^{1/k}}^
{x^{1/k}}\frac{dv}{1+\vert v\vert}\sum\limits_{r\sim
R}{\sum\limits_{\chi}}^{*}\int\limits_{v}^{T_{0}+v}\left\vert
L'\left(\frac{1}{2}+it,\chi
\right)\right\vert^{4}\,dt+T_{0}R^{2}L^{4}\\
&\ll & L^{5}\max_{\vert N\vert \leq
x^{1/k}}\int\limits_{N/2}^{N}\frac{dv}{1+\vert v\vert}\sum
\limits_{r\sim R}{\sum\limits_{\chi}}^{*}\int\limits_{v}^{T_{0}+v}
\left\vert L'\left(\frac{1}{2}+it,\chi\right)\right\vert^{4}\,dt+
T_{0}R^{2}L^{4}\\
&+&L^{5}\max_{\vert N\vert \leq
x^{1/k}}N^{-1}\int\limits_{0}^{T_{0}}dt \sum\limits_{r\sim
R}{\sum\limits_{\chi mod r}}^{*}
\int\limits_{\frac{N}{2}+t}^{N+t}\left\vert
L'\left(\frac{1}{2}+iv,\chi\right)\right\vert^{4}
\,dv\, +T_{0}R^{2}L^{4}\\
&\ll &  R^{2}T_{0}L^{c},\end{eqnarray*} Using lemma 3.8, (4.6) and
H\"older`s inequality we obtain
\begin{eqnarray*}& &\sum\limits_{r\sim R}{\sum\limits_{\chi}}^{*}
\int\limits_{0}^{T_{0}}\left\vert
F_{k}\left(\frac{1}{2}+it,\chi\right)
\right\vert\,dt\nonumber\\
&\ll& \left(\sum\limits_{r\sim R}{\sum\limits_{\chi}}^{*}
\int\limits_{0}^{T_{0}}\left\vert
f_{k,1}\left(\frac{1}{2}+it,\chi\right) \right\vert
dt\right)^{1/4} \left(\sum\limits_{r\sim
R}{\sum\limits_{\chi}}^{*} \int\limits_{0}^{T_{0}}\left\vert
f_{k,2}\left(\frac{1}{2}+it,\chi\right)
\right\vert dt\right)^{1/4}\\
&\times & \left(\sum\limits_{r\sim R}{\sum\limits_{\chi}}^{*}
\int\limits_{0}^{T_{0}}\left\vert \prod\limits_{j=3}^{10}
f_{k,j}\left(\frac{1}{2}+it,\chi\right)
\right\vert dt\right)^{1/2}\\
&\ll& (R^{2}T_{0})^{1/2}\left(
R^{2}T_{0}+\frac{x^{1/k}}{N_{k,1}N_{k,2}}\right)^{1/2} L^{c}\ll
x^{1/2k}R^{a_{k}}L^{-B},\end{eqnarray*}
by the definition of $T_{0}$ and the condition of the proposition.\\
{\bf Proposition \(2\)}\,\,\,{\it Let \(J=\{1,..,10\}\). If \(J\)
can be divided into two non overlapping subsets \(J_{1}\) and
$J_{2}$ such that\[\max \left (\prod\limits_{j\in
J_{1}}N_{k,j},\prod\limits_{j\in J_{2}}N_{k,j}\right)\ll
x^{\frac{1}{k}}P^{-2+2a_{k}-3\epsilon}\]
then (4.4) is true.}\\
{\it Proof:}\,\,\,Let \[F_{k,i}(s,\chi)=\prod\limits_{j\in
J_{i}}f_{k,j} (s,\chi)=\sum\limits_{n\ll M_{i}}b_{i}(n)
\chi(n)n^{-s},\enspace b_{i}(n)\ll d^{c}(n),\enspace i=1,2,\]
where $M_{i}=\prod\limits_{j\in J_{i}}N_{k,j},\,i=1,2.$ Applying
lemma 3.8, (4.3) and (4.6) we see \begin{eqnarray*}&
&\sum\limits_{r\sim
R}{\sum\limits_{\chi}}^{*}\int\limits_{0}^{T_{0}}\left\vert
F_{k}\left (\frac{1}{2}+it,\chi\right)\right\vert\,dt\\&\ll&
\left(\sum\limits_ {r\sim
R}{\sum\limits_{\chi}}^{*}\int\limits_{0}^{T_{0}}\left\vert
F_{k,1}\left(\frac{1}{2}+it,\chi\right)\right\vert\,dt\right)^{1/2}\left
({\sum\limits_{r \sim
R}\sum\limits_{\chi}}^{*}\int\limits_{0}^{T_{0}} \left\vert
F_{k,2}\left(\frac{1}{2}+it,\chi\right)\right\vert\,dt\right)
^{1/2}\\
&\ll &\left(R^{2}T_{0}+M_{1}\right)
^{1/2}\left(R^{2}T_{0}+M_{2}\right) ^{1/2} \nonumber\\& \ll &
R^{2}T_{0}+x^{\frac{1}{2k}}RP^{-1+a_{k}-\frac{3}{2}
\epsilon}T_{0}^{1/2}+x^{1/2k}L^{c}.\end{eqnarray*} This proves the
proposition because of $ R>L^{D}$.  Using proposition 1 and 2, we
can prove (\ref{eq:45}) in nearly the same way as (4.4) is proved
in \cite{b1}. The only difference in the proof is that instead of
assuming
\[ N_{k,i}N_{k,j}\leq P^{12/7+3\epsilon}\leq x^{2/5k},
\quad 1\leq i,j\leq 5,\,i\ne j\] as in \cite{b1}, we assume in
view of proposition 1 that
\[ N_{k,i}N_{k,j}\leq P^{2-2a_{k}+3\epsilon}\leq x^{2/5k},
\quad 1\leq i,j\leq 5,\,i\ne j.\] The proof of (\ref{eq:46}) goes
along the same lines. (4.1) is now proved in the case $R>L^{D}.$
The case $R\leq L^{D}$ is treated exactly as in \cite{b1}.
\section{Proof of lemma 3.6}\setcounter{equation}{0}
To prove the lemma it is enough to show that\[\max_{R\leq P/2}\sum
\limits_{r\sim
R}{\sum\limits_{\chi}}^{*}\max\limits_{\vert\lambda\vert\leq
1/rQ}\vert W_{k}(\l,\chi_{r}) \vert\ll
x^{1/k}R^{a_{k}}L^{A}.\]Arguing as in the section before - we do
not have to apply Gallagher`s lemma here -  we find
\begin{eqnarray*} W_{k}(\l,\chi)&\ll &L^{c}\max_{I_{{a_1}},..,I_{a_{2k+1}}}
\left\vert \int\limits_{-T}^{T}F(\frac{1}{2}+it,\chi)d\,t\int
\limits_{x/2^{k+1}}^{
x}u^{\frac{1}{2k}-1}e\left(\frac{t}{2k\pi}\log\,u+\l u\right)
d\,u\right\vert+x^{1/k}P^{-1},\end{eqnarray*} for $T=P^{3}.$
Estimating the inner integral by lemma 3.2 we obtain
\[\int\limits_{x/2^{k+1}}^{x} u^{\frac{1}{2k}-1}e\left(\frac{t}
{2k\pi}\log\,u+\l u\right)\,du ll x^{\frac{1}{2k}-1}
\min\left(\frac{x}{\sqrt{\vert t\vert +1}},\,
\frac{x}{\min\limits_{x/2^{k+1}<u\leq x}\vert t+2k\pi\l
u\vert}\right).\] Taking \(T_{0}=4k\pi x(rQ)^{-1}\) we conclude
that in order to prove the lemma it is enough to prove that for
$P\leq x^{\frac{7}{150}-\epsilon}$ and $2\leq k\leq 5$ there holds
\bea \sum\limits_{r\sim R}{\sum\limits_{\chi}}^{*}\int
\limits_{0}^{T_{0}}\left\vert
F_{k}(\frac{1}{2}+it,\chi)\right\vert\,dt &\ll&
x^{1/2k}R^{a_{k}}L^{c},
\\
\sum \limits_{r\sim R}{\sum\limits_{\chi}}^{*}
\int\limits_{T_{1}}^{2T_{1}} \left\vert
F_{k}(\frac{1}{2}+it,\chi)\right\vert \,dt &\ll& x^{1/2k}
R^{a_{k}}T_{1}L^{c},\quad T_{0}<\vert T_{1}\vert\leq T. \eea These
estimates are shown in the same way as (\ref{eq:45}) and
(\ref{eq:46}). Two
propositions analogous to the propositions 1 and 2 are proved:\\
{\bf Proposition \(3\)} \,\,\, {\it If there exist $N_{k,j_{1}}$
and $N_{k,j2}$ $(1\leq j_{1},j_{2}\leq 5)$ such that
$N_{k,j1}N_{k,j2}\geq
P^{2-2a_{k} +3\epsilon}$ then (5.1) is true.}\\
{\bf Proposition \(4\)}\,\,\,{\it Let \(J=\{1,..,10\}\). If \(J\)
can be divided into two non overlapping subsets \(J_{1}\) and
$J_{2}$ such that\[\max \left (\prod\limits_{j\in
J_{1}}N_{k,j},\prod\limits_{j\in J_{2}}N_{k,j}\right)\ll
x^{\frac{1}{k}}P^{-2+2a_{k}-3\epsilon}\]
then (5.1) is true.}\\
Remark: Here we do not need to treat the case $R>L^{D}$ separately
because we do not have to save a factor $L^{-B}.$


\begin{thebibliography}{99}
\bibitem{b} Bauer, C., {\it On a problem of the Goldbach-Waring Type}. Acta Mathematica Sinica, New
Series, vol. 14, No. 2 (1998), 223-234.
\bibitem{b1}Bauer, C. {\it An improvement on a theorem of the Goldbach-Waring
type,} Rocky Mountains Journal of Mathematics, Winter 2001, Volume
31, Number 4, 1151 - 1170.
\bibitem{g}Gallagher, P.X., {\it A large sieve density estimate near \(\sigma=1\).} Inventiones Math., 11 (1970), 329 - 339.
\bibitem{hb} Heath Brown, D.R., {\it Prime numbers in sort
intervals and a generalized Vaughan's identity.} Canadian J. Math., 34(1982),
1365 - 1377.
\bibitem{hua} Hua, L.K.; {\it Some results in the additive prime
number theory}, Quart. J. Math. (Oxford) 9(1938), 68 - 80.
\bibitem{lz_sq} Liu, J.Y.; Zhan, T.; {\it Squares of primes and
powers of two}, Monat. Math., 128(1999), 283 - 313.
\bibitem{mv}Montgomery, H.L.; Vaughan, R.C., {\it On the exceptional set in Goldbach`s problem.} Acta Arith., 27 (1975), 353 - 370.
\bibitem{pr}Prachar, K., {\it šber ein Problem vom Waring-Goldbach`schen Typ.} Monatshefte Mathematik., 57 (1953), 66 - 74.
\bibitem{t} Titchmarsh, E.C.; {\it The theory of the Riemann Zeta - Function.}
Second edition, Oxford. Clarendon Press, 1986.
\bibitem{v} Vinogradov, I.M., {\it Estimation of certain trigonometric sums with prime variables.}
Izv. Acada. Nauk SSSR . Ser. Mat., 3 no. 4 (1939), 371 - 398.
\bibitem{z} Zhan, T., {\it On the representation of a large odd integer as a sum
of three almost equal primes.} Acta Mathematica Sinica, 7(1991), No. 3, 259 - 272.
\end{thebibliography}
\end{document}